\documentclass[11pt,twoside]{amsart}

\usepackage[matrix,arrow]{xy}

\author{Florin Ambro} 
\address{Department of Mathematics\\
The Johns Hopkins University\\
3400 N. Charles, Baltimore MD 21218}
\email{ambro@chow.mat.jhu.edu}

\newcommand{\bP}{{\mathbb P}}
\newcommand{\bQ}{{\mathbb Q}}
\newcommand{\bZ}{{\mathbb Z}}
\newcommand{\bN}{{\mathbb N}}
\newcommand{\bR}{{\mathbb R}}

\newcommand{\bA}{{\mathbb A}}

\newcommand{\cF}{{\mathcal F}}
\newcommand{\cL}{{\mathcal L}}

\newcommand{\cO}{{\mathcal O}}

\newcommand{\cC}{{\mathcal C}}
\newcommand{\cI}{{\mathcal I}}

\newcommand{\cR}{{\mathcal R}}

\theoremstyle{plain}
\newtheorem{thm}{Theorem}[section]

\newtheorem{lem}[thm]{Lemma}
\newtheorem{cor}[thm]{Corollary}

\newtheorem{prop}[thm]{Proposition}

\theoremstyle{definition}
\newtheorem{defn}[thm]{Definition}
\newtheorem{defnprop}[thm]{Definition-Proposition}

\newtheorem{notation}[thm]{Notation}

\newtheorem{saynum}[thm]{}
\newtheorem{exmp}[thm]{Example}
\newtheorem{example}{Example}

\newtheorem{rem}[thm]{Remark}

\newtheorem{ack}{Acknowledgments}   

\theoremstyle{remark}

\newenvironment{sketch}{\begin{proof}[Sketch of proof]}{\end{proof}}

\begin{document}
\bibliographystyle{amsalpha+}
\title[The locus of log canonical singularities]
{The locus of log canonical singularities}
\maketitle


\tableofcontents

\setcounter{section}{-1}

\section{Introduction}

Given a log variety $(X,B_X)$, there is a naturally defined 
closed subscheme $LCS(X,B_X) \subset X$ \cite[3.14]{Sho}, called the 
{\em locus of log canonical singularities\/}.
The defining ideal sheaf, introduced by V. Shokurov (unpublished), 
is the algebraic counterpart of the multiplier ideal sheaves associated 
with singular hermitian metrics.
 
 Although $LCS(X,B_X)$ was initially introduced to measure how far
the log variety was from being log terminal (some authors called it
the non-Kawamata log terminal locus), it was realized recently
that $LCS(X,B_X)$ is worth studying in its own,
being an intermediate step for inductive arguments in higher dimensional
algebraic geometry. This technique was the main ingredient in the proof
of some of the basic theorems of the (Log) Minimal Model 
Program \cite[Ch. 2-4]{KMM}.

Our goal is to investigate the LCS locus in the ambient variety, 
and relate its singularities to those of the ambient space. 
As a first step in this direction, V. Shokurov proved \cite[3.6-8]{Sho} 
that $LCS(X,B_X)$ is normal if $(X,B_X)$ has pure log terminal 
singularities, and it has normal 
components intersecting normally (hence seminormal) if $(X,B_X)$ is 
strictly log terminal. The latter was generalized by J. Koll\'{a}r
\cite[17.5]{Kol1} to the log terminal case. 
\newline
In this paper, as conjectured by V. Shokurov, we show that 
$LCS(X,B_X)$ is seminormal if $(X,B_X)$ has log canonical singularities.

The natural category in which to study the LCS locus is that
of {\em (relative effective) log pairs\/} $\pi:(X,D)\to S$, where $D$ is 
{\em effective over $S$\/}, 
that is the negative part of $D$ is $\pi$-exceptional. 
A {\em log variety\/} is a log pair 
$(X,D)$ with $\pi=id_X$.
\\
\\
 We recall in Section 1 the basic definitions. In Section 2 we
introduce the LCS ideal sheaf for {\em log pairs\/}, which is a 
slight modification of V. Shokurov's definition. 
The ideal sheaves are isomorphic for log varieties
[Remark \ref{rem:ideals_coincide}]. The main technical
result is Theorem \ref{thm:ses}, which is a generalization of 
\cite[3.6]{Sho} (see also \cite[17.4]{Kol1}), and follows closely 
their proof. 
As a corollary, we obtain the contraction
which, via a formal seminormality result, implies the seminormality
of the LCS locus. 
In Section 3
we use Kawamata's technique to show that any finite union
of lc centers is seminormal, under some restrictions. Section 4
is an appendix on seminormality. 

\begin{ack}
  I would like to thank Professor V. Shokurov for setting up the 
  problem and also for his valuable support.
\end{ack}


\section{Log varieties and (relative effective) log pairs}
\label{sec:Log varieties and (relative effective) log pairs}

A variety is a reduced scheme of finite type over a fix field $k$.
We have to assume $char(k)=0$, since we use Kawamata-Viehweg 
vanishing as a main technical tool.

We first define the basic objects of this paper. 
\begin{defn}
   \begin{enumerate}
     \item A {\em relative effective log pair\/} $\pi:(X,D)\to S$
  is a normal variety $X$ equipped with an $\bR$-Weil divisor
  $D=\sum d_i D_i\ (d_i \in \bR)$, and a morphism $\pi:X \to S$ 
  such that
        \begin{enumerate}
    \item $K_X+D$ is an $\bR$-Cartier Weil divisor.
    \item $D$ is {\em relative effective\/}, that is the components
          $D_i$ of $D$ with negative coefficients are 
          $\pi$-exceptional ($codim(\pi(D_i),S)\ge 2$).
    \item $\pi$ is a {\em contraction\/}, that is $\cO_S=\pi_*\cO_X$. 
        \end{enumerate}
      \item A {\em log variety\/} is a log pair $\pi:(X,D)\to S$ such that
    $S=X,\pi=id_X$. Then the second condition is equivalent to $D$ being
    an effective divisor.
   \end{enumerate}
\end{defn}
We call $D$ the {\em pseudo-boundary\/} of the log pair, and we also
call $K+D$ a {\em log divisor\/}, since its sections correspond to 
rational differentials with poles along $D$. Recall that $K_X$ is a 
$\bZ$-Weil divisor on $X$, uniquely defined up to linear equivalence 
in its class, called the {\em canonical class\/}.
\newline
For simplicity of terminology, we will drop the adjective
``relative effective'' and if there is no danger
of confusion, we will also drop $\pi$ and $S$ from the notation, so
we will say that $(X,D)$ is a {\em log pair\/}.

\begin{saynum}
Given a log pair $\pi:(X,D)\to S$, a desingularization $\mu:Y \to X$
determines canonically a log pair $\varphi:(Y,D^Y)\to S$ \cite[pp.114]{Sho}.
\[ \xymatrix{
       (Y,D^Y) \ar[rr]^{\mu} \ar[dr]_{\varphi} &  & (X,D) \ar[dl]^{\pi} \\
         &                    S  &
} \]

Consider the following equality of Weil divisors:
$$
K_Y+\mu^{-1}D+\sum E_i=\mu^*(K_X+D)+\sum a_i E_i,
$$
where the sum runs over the $\mu$-exceptional divisor of $Y$, 
$\mu^{-1}D$ is the proper transform of the Weil divisor $D$.
The above formula determines uniquely the coefficients 
$a_i=a(E_i;X,D)$, which are called the {\em log discrepancy} of the 
exceptional divisors $E_i$. They are independent of the ambient resolution in 
which $E_i$ seats. We also extend the definition to non-exceptional
divisors $E$, declaring $a(E;X,D)=1-e$, where $e$ is the coefficient
of $E$ in $D$.  

If we denote $D^Y=\mu^{-1}D+\sum (1-a_i)E_i, \ \varphi=\pi\circ \mu$,
then $\varphi:(Y,D^Y) \to S$ becomes a log pair such that
$$K_Y+D^Y=\mu^*(K_X+D).$$
We say that $\mu$ is a {\em crepant\/} morphism of log pairs if the above 
equality holds.
Note that $(X,D)$ and $(Y,D^Y)$ have the same log discrepancies, and 
they should be viewed as being equivalent.
\end{saynum}

\begin{example}
Given a log variety $(X,B)$ and a resolution of singularities $\mu:Y \to X$,
it is easy to see that $\mu:(Y,B^Y) \to X$ is a log pair, while
$B^Y$ may have negative coefficients. This is the main example of log 
pairs appearing in the study of log varieties.
\end{example}

\begin{defn}
 A log pair $(X,D)$ is {\em log canonical\/} if all log
 discrepancies are nonnegative. In particular, $d_i \le 1 \ \forall i$
 (we say that $D$ is a {\em subboundary\/} in this case) .
\end{defn}

\begin{defn}
  \begin{enumerate}
    \item A log pair $(X,D)$ has {\em (log) nonsigular support\/}
  if $X$ is a nonsingular variety and $D=\sum d_i D_i$ is a Weil divisor
 such that $\cup_{d_i \ne 0} D_i$ is a union of smooth divisors
 intersecting transversely.
     \item A {\em log resolution \/} of a log pair $(X,D)$
 is a log pair induced on a resolution of singularities $(Y,D^Y)$
 which has nonsingular support.
  \end{enumerate}
\end{defn}
\begin{exmp}
 Assume $(X,D)$ is a log pair with nonsingular support.
 Then $(X,D)$ is log canonical iff 
 $d_i \le 1, \ \forall i$ (see the proof of \ref{defn-prop:LCS}).
\end{exmp}

\begin{exmp}
 Let $X$ be a toric variety and $B_X=\sum B_i$ be the complement
 of the embedded torus. Then $(X,B_X)$ is a Calabi-Yau log variety:
 $$K_X+B_X=0.$$
 Moreover, $(X,B_X)$ is log canonical 
 (see \cite[4.8]{YPG}, or \cite[3.1]{Ale1}).
\end{exmp}

\begin{rem}
 Although log canonicity involves all possible prime divisors 
 with center on $X$, it is enough to check it on a {\em log resolution\/}
 $\mu:(Y,D^Y) \to (X,D)$. 
 From definition, $(X,D)$ is log canonical iff $(Y,D^Y)$ 
 is log canonical which, in turn, is equivalent to $D^Y$ being a subboundary. 
\end{rem}


\section{The locus of log canonical singularities}
\label{sec:The locus of log canonical singularities}

\begin{notation} For an $\bR$-Weil divisor  $D=\sum d_i D_i, d_i \in \bR$
on a normal variety $X$ we define
 \begin{enumerate}
   \item the coherent {\em divisorial sheaf\/} 
       ${\cO}_X(D) \subset {\cR}_X =K(X)$ defined as
$$
H^0(U,{\cO}_X(D))=\{f \in K(X); (f)+D|_U\ge 0\},\ \ U \subseteq X.
$$
   \item If $\cF$ a coherent sheaf on $X$, and ${\cO}_X(D)$ is an 
  invertible sheaf, we denote $\cF(D):=\cF \otimes \cO_X(D)$.
   \item the {\em round up (down) of\/} $D$, 
          $\lceil D \rceil=\sum \lceil d_i \rceil D_i$ 
           ($\lfloor D \rfloor=\sum \lfloor d_i \rfloor D_i$).
   \item the {\em positive (negative) part of\/} $D$, 
            $D^+=\sum_{d_i>0}d_i D_i$ ($D^-=\sum_{d_i>0}d_i D_i$), so
         the decomposition $D=D^+ + D^-$ holds.
 \end{enumerate}
Note the identities 
 \begin{itemize}
   \item[a)] ${\cO}_X(D)={\cO}_X(\lfloor D\rfloor)$.
   \item[b)] $\lceil -D \rceil=-\lfloor D\rfloor$. 
   \item[c)] $\lceil-(D^+) \rceil=-\sum_{d_i \ge 1} 
    \lfloor d_i \rfloor D_i$.
 \end{itemize}
We declare that taking the positive (negative) part of a divisor
has precedence over all other operations. For example, we will write
$-D^-$ for $-(D^-)$.
\end{notation}
 
\begin{defnprop}
 \label{defn-prop:LCS} 
Let $(X,D)$ be a log pair and let 
$\mu:Y \to X$ a log resolution with $D^Y$ the corresponding
pseudo-boundary on $Y$. Then the coherent ideal sheaf on $X$  
$$
{\cI}(X,D)=\mu_* {\mathcal Hom}_Y({\cO}_Y((D^Y)^+),{\cO}_Y)
=\mu_* {\cO}_Y (\lceil-(D^Y)^+ \rceil)
$$
is independent of the log resolution. 
The induced subscheme of $X$, denoted $LCS(X,D)$, is called {\em the
locus of log canonical singularities of the log pair\/} $(X,D)$. 
\end{defnprop}
\begin{rem} The above definition is a slight modification
of V. Shokurov's definition of the LCS ideal (see also
\cite[2.16]{Kol2}). He defined the coherent sheaf 
${\cI'}(X,D)=\mu_* {\cO}_Y (\lceil -D^Y \rceil)$,
which is isomorphic to an ideal sheaf only if $D$ is effective,
in which case ${\cI'}(X,D) \simeq {\cI}(X,D)$ 
(see Remark \ref{rem:ideals_coincide}). 
It is interesting that ${\cI'}(X,D)$
has good vanishing properties by the very definition, which is
not the case for the actual ideal ${\cI}(X,D)$.
\end{rem}
\begin{proof}
Using Hironaka's hut, it is enough to check that $(X,D)$ has
nonsingular support, and $\tau:Y \to X$ is a sequence of blow-ups 
with nonsingular centers, then
$$
\tau_*{\cO}_Y(\lceil -(D^Y)^+ \rceil)={\cO}_X(\lceil -D^+ \rceil).
$$
Indeed, let $\{E_i\}_{i=1}^t$ be the exceptional locus of $\tau$, with
$m_i=cod_X(\tau(E_i))$. Then $K_Y=\tau^*K_X+ \sum_{i=1}^t \alpha_i E_i$
and $D^Y= \tau^*D- \sum_{i=1}^t \alpha_i E_i$, where $\alpha_i \ge m_i-1$. 
We claim that 
$
\lfloor D^Y \rfloor
\le \tau^*\lfloor D \rfloor
$.
\newline
Indeed, if $\tau(E_k)$ lies in 
$D_1,\ldots,D_s \ (s \le t)$ only, then 
$m_k \ge s$ and $\lfloor \sum_{j=1}^s d_j \rfloor \le
\sum_{j=1}^s \lfloor d_j \rfloor +s-1 \le 
\sum_{j=1}^s \lfloor d_j \rfloor + \alpha_k$, that is
$
\lfloor \tau^*D \rfloor \le
\tau^*\lfloor D \rfloor + \sum_{i=1}^t \alpha_i E_i
$.
\newline
Now $\lfloor (D^Y)^+ \rfloor={\lfloor D^Y \rfloor}^+ \le
(\tau^*\lfloor D \rfloor)^+ \le \tau^*( {\lfloor D \rfloor}^+)=
  \tau^*( \lfloor D^+ \rfloor)$, hence
$$
\lfloor (D^Y)^+ \rfloor \le
\tau^*\lfloor D^+ \rfloor.
$$ 
This last inequality, together with the fact that $D$ has simple
normal crossings support implies at once that
$\tau_*{\cO}_Y(-\lfloor (D^Y)^+ \rfloor)={\cO}_X(-\lfloor D^+ \rfloor)$.
\end{proof}

We have a dichotomy: either $LCS(X,D)=\emptyset$,
in which case we say that $(X,D)$ has {\em Kawamata log terminal
singularities\/} ({\em klt\/} for short), or 
$LCS(X,D)\ne \emptyset$ is a proper subscheme of $X$. We are mainly 
interested in the second case. 
\begin{example}
 Let $(X,D)$ be a log pair with log nonsingular
 support. If we write 
  $D=\sum d_i D_i$, and $E= \sum_{d_i \ge 1} \lfloor d_i \rfloor D_i$,
then
$$LCS(X,D)=(E,\cO_E).
$$
\end{example}

\begin{example}
Let $S=\bA^2$ and let $C:(y^2-x^3=0) \subset S$ a curve with a cusp
at the origin $P$. Then 
 \begin{enumerate}
   \item $LCS(S,t C)=\emptyset$ if $0 \le t <\frac{5}{6}$.
   \item $LCS(S, \frac{5}{6} C)=\{P\}$ and $(S,\frac{5}{6} C)$ 
               is log canonical.
   \item $LCS(S,t C)=C$ as a set and $(S,t C)$ is not log canonical
               for $t>\frac{5}{6}$.
  \end{enumerate}
\end{example}
\begin{example} Let $L,H \subset \bP^3$ be a line and a plane intersecting
in a point. Let $H_1,H_2,H_3$ be three general planes passing through
the line $L$. Let $B=H+\frac{2}{3}(H_1+H_2+H_3)$. Then
$(\bP^3,B)$ is a log variety with log canonical singularities, and
$$LCS(\bP^3,B)=L \cup H.$$
This is an example of non pure dimensional LCS locus.
\end{example}

\begin{example} Let $L_1,L_2,L_3 \subset \bP^3$ be three 
lines passing through a point $P$ such that 
$\dim T_P(L_1 \cup L_2 \cup L_3)=3$. Let $H_{i1},H_{i2},H_{i3}$
be generic planes passing through the line $L_i$, for $1\le  i \le 3$.
Let $B=\frac{2}{3}\sum H_{ij}$. Then $(\bP^3,B)$ has log canonical 
singularities, and
$$LCS(\bP^3,B)=L_1 \cup L_2 \cup L_3.$$
The same is true if we consider $n$ general lines passing through a point
in $\bP^n$.
\end{example}
\begin{rem} Let us denote $X^{(j)}=\{\eta \in X; codim(\bar{\eta},X)=j\}$.
Assume the log pair $(X,D)$ is log canonical. 
It is interesting that the $LCS(X,D)$ 
is smooth in $X^{(1)}$, it has at most ordinary
double points in $X^{(2)}$, and it
can have triple ordinary points in $X^{(3)}$
(similarly for any $n>3$). So we could naturally ask if
the triple ordinary points are the only type of singularities appearing
in $X^{(3)}$. Although they are seminormal, the answer might
be negative!
\end{rem}

\begin{prop} Let $\tau:(Y,D_Y) \to (X,D_X)$ be a birational contraction
of log pairs and assume that $\tau$ is crepant, that is
$$
\tau^*(K_X+D_X)=K_Y+D_Y.
$$
Then $\tau$ induces a dominant morphism between the LCS schemes
$$
\tau':LCS(Y,D_Y) \to LCS(X,D_X),
$$
\end{prop}
\begin{proof} Let $\mu:Z \to Y$ be a log resolution. 
Then $\tau \circ \mu:Z \to X$ is also a log resolution, and 
$(D_X)^Z=(D_Y)^Z$. Therefore ${\cI}(X,D_X)=\tau_* {\cI}(Y,D_Y)$. This
easily implies that ${\cO}_{ LCS(X,D_X)} \to {\cO}_{ LCS(Y,D_Y)}$
is injective, hence $\tau'$ is dominant. 
\end{proof}

\begin{thm}
 \label{thm:ses}
 Let $\pi:(X,D)\to S$ be a log pair such that
 $-(K_X+D)$ is $\pi$-nef and $\pi$-big. 
 \begin{enumerate}
   \item  The following sequence is exact:
 \[ \xymatrix{
    0 \ar[r] & \pi_*{\cI}(X,D) \ar[r]  &
     {\cO}_S \ar[r] & \pi_*{\cO}_{LCS(X,D)} \ar[r] & 0
 } \]
   \item $R^1{\pi}_*{\cI}(X,D)=0$ if $\pi$ is a birational contraction,
          or if $R^1 \pi_*\cO_X=0$. 
  \end{enumerate} 
 \end{thm}

\begin{proof} We assume first that $(X,D)$ has log nonsingular
 support.
 Denoting ${\cI}={\cI}(X,D)$ and $E=\lceil -D^- \rceil$ we have
 $$
 {\cO}_X(\lceil -D \rceil)={\cI} \otimes {\cO}_X(E).
 $$
 Since $-(K_X+D)$ is $\pi$-nef and $\pi$-big, and $(X,D)$ is log nonsingular,
 Kawamata-Viehweg vanishing implies 
 $R^j \pi_*{\cO}_X(\lceil -D \rceil)=0,\ \forall j \ge 1$, that
 is 
 $$
 R^j \pi_*{\cI}(E)=0, \ \forall j \ge 1.
 $$
 Look at the following commutative diagram with exact rows:
\[ \xymatrix{
   0 \ar[r] & {\pi_*{\cI}} \ar[d]_{j_0} \ar[r] &
    {\pi_*{\cO}_X} \ar[d]_{j_1} \ar[r]^{i_1} & 
   {\pi_*{\cO}_{LCS}} \ar[d]_{j_2} \\
   0 \ar[r] & \pi_*{\cI}(E) \ar[r] &
   \pi_*{\cO}_X(E) \ar[r]^{i_2} &  \pi_*{\cO}_{LCS}(E) \ar[r] &
   R^1\pi_*{\cI}(E)=0
} \]
Since $E$ is effective $\pi$-exceptional, $j_1$ is an isomorphism. Moreover,
$i_2$ is surjective due to vanishing, hence $j_2$ is surjective.
But $j_2$ is injective, hence $j_0,j_1,j_2$ are all isomorphisms
and $i_1$ is surjective. This proves the first part. 
\newline
For the second, let us assume that $\pi$ is birational. 
Consider now the following comutative diagram
\[ \xymatrix{
     0 \ar[r] & R^1\pi_*{\cI} \ar[d] \ar[r] & R^1\pi_*{\cO}_X \ar[d]^j\\
              & R^1\pi_*{\cI}(E)=0 \ar[r] & R^1\pi_*{\cO}_X(E)
} \]
where the top row is exact from the above argument.
We claim that $j$ is injective. Indeed, $\pi_*{\cO}_E(E)$ surjects onto 
$Ker(j)$ and $\pi_*{\cO}_E(E)=0$ \cite[1-3-2]{KMM}. Therefore
the morphism $R^1\pi_*{\cI} \to R^1\pi_*{\cO}_Y(E)$ is 
injective too, hence $R^1\pi_*{\cI}=0$.
\\
\\
 Now, for the general case, let $\mu:(Y,D^Y) \to (X,D)$ be a log resolution
and denote $\nu=\pi \circ \mu$ and $E=\lceil -(D^Y)^- \rceil$. 
We have the following diagram with exact rows:
\[ \xymatrix{
   0 \ar[r] & \pi_*{\cI}(X,D) \ar[d]_{=} \ar[r] &
    {\pi_*{\cO}_X} \ar[d]_{\simeq} \ar[r] & 
   \pi_*{\cO}_{LCS(X,D)} \ar[d] \\
   0 \ar[r] & \nu_*{\cI}(Y,D^Y) \ar[r] &
   \nu_*{\cO}_Y \ar[r] &  \nu_*{\cO}_{LCS(Y,D^Y)} \ar[r] & 0
} \]
where the bottom row is exact from the previous step. With
the same argument as above, we obtain that all the vertical arrows
are isomorphisms. In particular, the last arrow of the top arrow is
surjective. Finally, the exact sequence of lower terms of the
spectral sequence
$$
E^{p,q}_2=R^p \pi_* R^q \mu_* \cI(Y,D^Y) \Longrightarrow 
         E^{p+q}=R^{p+q}\nu_*\cI(Y,D^Y)
$$
gives the injection $R^1\pi_* \cI(X,D) \hookrightarrow R^1\nu_*\cI(Y,D^Y)$, 
hence we are done from the previous case.

\end{proof}
\begin{rem} 
  \label{rem:ideals_coincide}
Note that if $D$ is effective, $j_0$ is an isomorphism 
between ${\cI}(X,D)$ and ${\cI'}(X,D)$.
\end{rem}

\begin{prop}  

 Let $\tau:(Y,D_Y) \to (X,D_X)$ be a crepant birational contraction
 of log pairs
 \[ \xymatrix{
       (Y,D_Y) \ar[rr]^{\tau} \ar[dr]_{\varphi} &  & (X,D_X) \ar[dl]^{\pi} \\
         &                    S  &
 } \]
 Then $\varphi_*{\cO}_{LCS(Y,D_Y)}=\pi_*{\cO}_{LCS(X,D_X)}$, that is
 $\pi_*{\cO}_{LCS(X,D_X)}$ is a birational invariant of log pairs.
 \newline
 In particular, if $(X,B)$ is a log variety and $\tau:Y \to X$
 is a resolution, then $\tau': LCS(Y,B^Y) \to LCS(X,B)$ 
 is a contraction, that is
 $
 {\cO}_{LCS(X,B)}=\tau'_*{\cO}_{LCS(Y,B^Y)}
 $.
 \end{prop}

 \begin{proof}Let $\mu:Z \to Y$ be a log resolution. 
Then $\tau \circ \mu:Z \to X$ is also a log resolution, and 
$(D_X)^Z=(D_Y)^Z$. We then apply the previous theorem for
$\tau \circ \mu:Z \to X$.
 \end{proof}

\begin{cor}[Connectedness Lemma {\cite[5.7]{Sho}},{\cite[17.4]{Kol1}}] 
 Assume \newline $\pi:(X,D)\to S$ is a log pair such that 
 $-(K_X+D)$ is $\pi$-nef and $\pi$-big.
 Then $LCS(X,D) \cap \pi^{-1}(s)$ is connected for every $s \in S$.
\end{cor}

\begin{proof}
 The surjection 
 $$
 {\cO}_S \to \pi_*{\cO}_{LCS(X,D)} \to 0.
 $$ 
 easily implies the connectivity of the fibers.
\end{proof}

\begin{lem}
  \label{lem:red_lcs_issn}
 Let $(X,B)$ be a log variety.
 Assume there is a log resolution $\mu:(Y,B^Y) \to (X,B)$ such that
 $LCS(Y,B^Y)$ is a reduced scheme. Then $LCS(X,B)$ is seminormal. 
\end{lem}
\begin{proof}
Note that $LCS(Y,B^Y)$ is reduced iff  $\lfloor (B^Y)^+ \rfloor$ is a reduced
divisor. Then $LCS(Y,B^Y)$
is a simple normal crossings divisor with the induced reduced 
structure, which is seminormal by \ref{lem:snc_is_sn}. But $LCS(X,B)$
is a contraction of $LCS(Y,B^Y)$, so we can apply \ref{prop:contr_preserve_sn}.
\end{proof}
\begin{cor}
Assume $(X,B)$ is a log variety with log canonical singularities.
Then $LCS(X,B)$ is a seminormal variety.
\end{cor}

\begin{thm}[V. Shokurov]
Let $(X,B)$ be a log variety and assume that $B$ is effective.
Let $\pi:X \to S$ be a proper morphism and let $L$ be a
Cartier divisor on $X$ such that
$$
L \equiv K_X+B+H
$$
where $H$ is a $\pi$-nef and $\pi$-big $\bR$-Cartier divisor. Then
$$
R^j \pi_* ({\cI}(X,B)(L))=0,\ \forall j\ge1.
$$
\end{thm}

\begin{proof}
Let $\mu:(Y,B^Y) \to (X,B)$ be a log resolution and denote 
$\nu=\pi \circ \mu$. We can assume that there is an effective 
$\bR$-divisor $F$, and a $\nu$-ample $\bR$-divisor $A$
on $Y$ such that $\mu^*H=F+A$ and $\lceil -B^Y-F \rceil=\lceil -B^Y\rceil$.
Since $-(K_Y+B^Y)$ is $\mu$-nef and $\mu$-big, Kawamata-Viehweg vanishing
gives $R^j \mu_* {\cO}_Y(\lceil -B^Y\rceil)=0,\ \forall j\ge1$, hence,
by the projection formula,
$$
R^j \mu_* {\cO}_Y(\lceil -B^Y\rceil +\mu^*L)=0,\ \forall j\ge1.
$$
Therefore the Grothendieck spectral sequence degenerates and 
the following isomorphism holds
$$
R^j \pi_* (\mu_*{\cO}_Y(\lceil -B^Y\rceil)(L))
\simeq 
R^j \nu_* {\cO}_Y(\lceil -B^Y\rceil + \mu^*L),\ \forall j\ge 0.
$$
But $\mu^*{\cL}-K_Y-B^Y-F$ is $\nu$-ample, so Kawamata-Viehweg vanishing
gives
$$
R^j \nu_* {\cO}_Y(\lceil -B^Y\rceil + \mu^*L)=0,\ \forall j\ge1.
$$
Finally, since $B$ is effective,
$\mu_*{\cO}_Y(\lceil -B^Y\rceil)={\cI'}(X,D) \simeq {\cI}(X,D)$,
so we obtain our vanishing.
\end{proof}



\section{Kawamata's lc centers and perturbation trick}
\label{sec:Kawamata's lc centers and perturbation trick}

Kawamata's pertubation technique applies for log varieties
$(X,B)$ with the following property \cite{Ka2}:
there is another log variety $(X,B^o)$ such that
 \begin{enumerate}
   \item $(X,B^o)$ has Kawamata log terminal singularities,
   \item $B^o < B$.
 \end{enumerate}
We will assume this throughout this section.

\begin{defn}[\cite{Ka1}] Let $(X,B)$ be a log variety. Then
any codimension one component $E$ of $LCS(Y,B^Y)$, for
any resolution $\mu:Y \to X$, is called a {\em log canonical (lc) place\/}.
The image on $X$ of a log canonical place is 
called a {\em log canonical (lc) center\/}.
\end{defn}
Log canonical centers are building blocks of the LCS locus.
Note that the LCS locus is the union of all the lc centers. All
the irreducible components of the LCS locus are lc centers.

\begin{prop} 
 \label{prop:union_centers_issn}
Assume $(X,B)$ is a log canonical variety.
Then any finite union of lc centers is a seminormal variety. 
In particular, every irreducible component of $LCS(X,B)$
is seminormal.
\end{prop} 

\begin{proof}(cf. \cite[1.5]{Ka1})
 Let $W_1,\ldots,W_k$ be lc centers for $(X,B)$, 
 $W=\cup_i W_i$, and let $E_{ij}$
 be all corresponding lc places on a log resolution $\mu:Y \to X$
 such that $\mu(E_{ij})=W_i$. 
 Let $H_i\supseteq W_i (1 \le i \le k)$ be generic effective divisors.
 Define
 $$
 B_{\epsilon}=(1-\epsilon)B+\epsilon B^o+\sum \epsilon a_i H_i,\ \ 
 0 < \epsilon \ll 1, a_i \in \bR.
 $$
 Then $K_X+B_{\epsilon}=K_X+B+\epsilon(\sum a_i H_i -(B-B^o))$
 is an $\bR$-Cartier divisor. 
 Let $a_i\ (1 \le i \le k)$ be the smallest positive numbers satisfying the
 inequalities
 $$
 a(E_{ij};X,B_{\epsilon})\le a(E_{ij};X,B),\  \forall i,j, \ \forall \epsilon.
 $$ 
 ($a_i=min_j \frac{\nu(E_{ij};B-B^o)}{\nu(E_{ij};H_i)}$, where
 $\nu(E;M)$ denotes the coefficient of $E\subset Y$ in $\mu^*(M)$) 
 \newline
 Then $LCS(X,B_{\epsilon})= W$.
 Indeed, $LCS(X,B_{\epsilon})\subseteq LCS(X,B)$ for small $\epsilon$, and
 if $E \in LCS(Y,B^Y), \mu(E) \not\subseteq W$, then
 $a(E;X,B_{\epsilon})=a(E;B-\epsilon(B-B^o))<a(E;B)\le 0$.
 \newline
 Therefore $(X,B_{\epsilon})$ is log canonical in the generic points
 of $W=LCS(X,B_{\epsilon})$, and even if it has worse  
 singularities in proper points of $W$, $\lfloor B_{\epsilon}^Y \rfloor$ 
 is a reduced divisor on $Y$ for $\epsilon$ small enough. Therefore
 \ref{lem:red_lcs_issn} gives the seminormality of $W$.
 \end{proof}

\begin{prop}\cite[1.5]{Ka1}
 Assume the log variety $(X,B)$ has log canonical singularities. 
 Let $W_1,W_2$ be two lc centers of $(X,B)$ on $X$. 
 Then every irreducible component of $W_1 \cap W_2$ is a lc center
 for $(X,B)$. In particular, there are minimal (with respect to inclusion)
 centers.
\end{prop}

\begin{lem}[\cite{Ka1}] Let $(X,B)$ be a log canonical variety, and let
$W \subseteq LCS(X,B)$ be a minimal lc center. Then there
is a log canonical variety  $(X,B')$ and a log resolution
$\mu:(Y,{B'}^Y) \to (X,B')$ such that $E= LCS(Y,{B'}^Y)$ is a 
smooth prime divisor.
In particular, there is an induced contraction
$\nu:E \to W$, hence $W$ is normal.
\end{lem}

\begin{proof}
We can assume $W=LCS(X,B)$ using the argument of 
\ref{prop:union_centers_issn}. 
The variety $(X,B)$ stays log canonical because $W$ is minimal.
Let $\mu:(Y,B^Y) \to (X,B)$ be a log resolution, and let
$LCS(Y,B^Y)=\cup_{1\le i \le k} E_i$, where $\mu(E_i)=W$. 
\newline
We have to decrease the log dicrepancies of all but one of the $E_i$'s and
the following trick was kindly suggested by V. Shokurov.
Let $\{\{m_i\}_i\}$ be a bounded family of integer vectors with 
integers entries, and let $A$ be an ample Cartier divisor on $Y$ and
$H$ a nef and big Cartier divisor on $X$.
 After taking a high multiple of $A$, we can assume that $|A-\sum m_i E_i|$
is a free linear system for all vectors $\{m_i\}_i$ in our bounded family.
Since $\mu^*H$ is nef and big on $Y$, there is $N \in \bN$ such that
$\mu^*(NH)=E+A$, with $E$ effective. Therefore after scaling the family 
$\{\{m_i\}_i\}$ and $E$ with $N$, there is an effective $\bQ$-Cartier divisor 
$M \sim_{\bQ} H$ such that
$$
\mu^*M = E + \sum m_i E_i + F, \forall \{m_i\}_i,
$$
where $E$ is an effective $\bQ$ -divisor (same for all $\{m_i\}_i$), and $F$ 
is a $\bQ$-free effective divisor, not containing $E_i$'s in its support.
Define
$$
B_{\epsilon}=(1-\epsilon)B+\epsilon B^o+ \epsilon aM, \ 
0 < \epsilon \ll 1, a \in \bR.
$$ 
Then $B_{\epsilon}^Y=B^Y+ \epsilon(aE + \sum  a m_i E_i - \mu^*(B-B^o))+ 
\epsilon a F$.
Let $a$ be the smallest positive number satisfying the inequalities 
$$
a(E_i;X,B_{\epsilon})\le a(E_i;X,B),\  \forall i, \ \forall \epsilon,
$$
that is $a=min_i \frac{\nu(E_i,B-B^o)}{m_i+\nu(E_i,E)}$.
Since we have a family, we can asssume that the 
equality holds for exactly one $E_i$. We just take now $B'=B_{\epsilon}$.
\end{proof}

\section{Appendix on seminormal varieties}
\label{sec:Appendix on seminormal varieties}
We say that a morphism $f:Y \to X$ is
a {\em quasi-isomorphism\/} if it is a universal homeomorphism such
that $k(f(x))\stackrel{\simeq}{\to} k(x)$ for all Grothendieck points 
$x \in X$. Note that any quasi-isomorphism is birational.
A {\em universal homeomorphism\/} is a morphism $f:Y \to X$ such that
for any base change $X' \to X$, the induced morphism $Y \times_{X'} X \to X'$
is a (topologically) homeomorphism.
\begin{defn}
\begin{enumerate}
    \item Let $f:Y \to X$ be a dominant morphism of preschemes such
that $f_*{\cO}_Y$ is a quasi-coherent ${\cO}_X$-algebra 
(this always happens in applications, for example if
$f$ is quasi-compact and quasi-separated). 
\newline
The {\em seminormalization of $f$\/} is an integral quasi-isomorhism 
$sn_f:X^{sn,f} \to X$ which factors $f$, and is maximal 
with respect to this property. That is, if $g:Z \to X$ is another 
integral quasi-isomorhism factoring $f$, then there is a
unique morphism $X^{sn,f} \to Z$ making the following diagram commutative:
\[ \xymatrix{
   Y \ar[dd]_f \ar[dr] \ar[drr] \\
   & X^{sn,f} \ar@{{}-->}[r] \ar[dl]_{sn_f} & Z \ar[dll]^g \\
   X
} \]
   \item We say that {\em $f$ is seminormal\/}, or
that $X$ is seminormal in $Y$,  if $sn_f$ is an isomorphism.
\end{enumerate}
\end{defn}

\begin{sketch}[See \cite{AB,Tr}] We have to prove existence only.
Let ${\cC}$ be the integral closure of ${\cO}_X$ in 
$f_*{\cO}_Y$, which is a quasi-coherent ${\cO}_X$-algebra (EGA II.6.3.4). 
Define ${\cC^{sn}}$ as follows
$$
H^0(U,{\cC^{sn}})=\{s \in H^0(U,{\cC});
   \forall x \in U, s_x \in {\cO}_x+R({\cC}_x) \},
 \ U \subseteq X$$
where $R(A)$ denotes the radical of the ring $A$, that is the intersection
of all its maximal ideals.
\newline 
Then it follows that ${\cC^{sn}}$ is a quasi-coherent ${\cO}_X$-algebra
and 
$$
sn_f:{\mathcal Spec}_X(\cC^{sn}) \to X
$$ 
satisfies the required universal property. For a detailed proof
with ``quasi-isomorphism'' replaced by ``universal homeomorphism'', 
see \cite{AB}.
\end{sketch}

\begin{saynum} (Functoriality)
Let $f_1:Y_1 \to X_1$ and $f_2:Y_2 \to X_2$ be
two dominant morphisms such that there are two morphisms 
$\alpha:X_1 \to X_2,\beta:Y_1 \to Y_2$ with 
$f_2\circ \beta=\alpha \circ f_1$. Then there is a unique morphism
$\alpha^*:X_1^{sn,f_1} \to X_2^{{sn,f_2}}$ such that 
the following diagram is commutative:
\[ \xymatrix{
   Y_1 \ar[dd] \ar[dr] \ar[rr] & & Y_2 \ar[dd] \ar[dr] &\\
   & X_1^{sn,f_1} \ar[dl]^{sn_{f_1}} \ar@{{}-->}[rr] & & 
     X_2^{sn,f_2} \ar[dl]^{sn_{f_2}} \\
   X_1 \ar[rr] & & X_2 &
} \]
Indeed, giving $\alpha^*$ is the same as giving a map $X_1^{sn,f_1} \to
Z=X^{sn,f_1}\times_{X_2} X_2^{sn,f_2}$. 
But $Z \to X_1$ is an integral quasi-isomorphism factoring $f_1$, 
hence the existence and uniqueness
of $\alpha^*$ follows from the universal property of $sn_{f_1}$.
\newline Using the above functoriality and chasing diagrams, it is easy
to see that seminormal morphisms behave well under composition.
If $f:Z \to Y$ and $g:Y \to X$ are dominant morphisms, then
$g,f \ sn \Longrightarrow g \circ f \ sn$, and 
$g \circ f \ sn \Longrightarrow g \ sn$.
\end{saynum}

\begin{saynum}(Contractions) We say that $f:Y \to X$ is a {\em contraction}
({\em fiber space\/}) if the natural morphism 
${\cO}_X \to f_*{\cO}_Y$ is an isomorphism (is an algebraically closed 
extension). Contractions and fiber spaces are example of
seminormal morphisms. Indeed, ${\cC}^{sn}={\cC}={\cO}_X$ in this case.
\end{saynum}

\begin{saynum}Let's fix a field $k$ of any characteristic. From now on we 
consider {\em varieties\/} only, i.e. algebraic reduced $k$-schemes
(possibly with more than one irreducible component).
Let $X$ be a variety with normalization $\pi:\bar{X} \to X$,
which is a birational finite morphism. Define 
{\sl the seminormalization of $X$} to be the seminormalization
of $X$ in $\bar{X}$.
Then $sn_X:X^{sn} \to X$ is a {\em finite 
 quasi-isomorphism\/} which is maximal in the following sense:
\begin{quote}
For any quasi-isomorphism $Z \stackrel{g}{\to} X$
from a variety $Z$, there is a unique morphism 
$\sigma:X^{sn}\to Z$ such that $g \circ \sigma=sn_X$.
\end{quote}

There is a functor associating to any variety $X$ its seminormalization
$X^{sn}$, and to any morphism $f:X \to Y$ its unique extension 
$f^{sn}:X^{sn} \to Y^{sn}$.
\[ \xymatrix{
   X^{sn} \ar[d] \ar[r]^{f^{sn}} & Y^{sn} \ar[d]\\
   X \ar[r]^f & Y
} \]
This follows from the general functoriality, since any morphism
lifts to normalizations.

\end{saynum}

\begin{prop}
 \label{prop:contr_preserve_sn}
Let $f:Y \to X$ be a contraction, or more generally, a seminormal 
morphism. Then $Y$ seminormal implies that $X$ is seminormal.
\end{prop}
\begin{proof}
Indeed, $sn_Y$ is an isomorphism and
$f=sn_X \circ (f^{sn}\circ sn_Y^{-1})$, so $f$ factors through the
quasi-isomorphism $sn_X$. But the seminormalization of
$X$ in $Y$ is $X$, hence $sn_X$ is an isomorphism.
\end{proof}

\begin{lem}
 \label{lem:snc_is_sn}
Let $D$ be the support of a reduced normal crossing divisor
on a nonsingular variety $X$. Then $D$ is seminormal.
\end{lem}
\begin{proof} Since a local ring ${\cO}$ is seminormal iff its 
completion ${\cO}^{-}$ is seminormal \cite[5.3]{GT}, we can assume 
${\cO}_{D,P}=k[X_1,\ldots ,X_n]/(X_1 \cdots X_s), \ \ s \le n$.
It is easy to see that this ring is seminormal.
\end{proof}

\begin{rem} By Serre's criterion, a variety $X$ is normal iff
  \begin{enumerate}
   \item $X$ is nonsingular in codimension $1$ and
   \item $X$ is $S_2$-saturated, that is
 $\cO_X = j_*\cO_{X-Z}$ for every closed subset $Z$ of $X$
 codimension at least $2$. 
   \end{enumerate}
Similarly, an $S_2$-saturated variety $X$ is seminormal iff $X$ is
seminormal in codimension $1$ \cite[2.6]{GT}. Moreover, the codimension $1$ 
seminormal singularities are classified. They basically look like
the origin on the $n$ coordinate axes in $\bA^n$. 
\end{rem} 


\ifx\undefined\bysame
\newcommand{\bysame}{\leavevmode\hbox to3em{\hrulefill}\,}
\fi

\end{document}